\documentclass[preprint,10pt]{elsarticle}
\pdfoutput=1
\pdfcatalog{/PageLayout /SinglePage}



\makeatletter
\def\ps@pprintTitle{%
   \let\@oddhead\@empty
   \let\@evenhead\@empty
   \def\@oddfoot{\reset@font\hfil\thepage\hfil}
   \let\@evenfoot\@oddfoot
}
\makeatother

\usepackage{titlesec}
\usepackage{xcolor}
\usepackage{enumitem}
\usepackage{hyperref}
\hypersetup{colorlinks=true,urlcolor=black,pdfpagemode=UseNone,pdfstartview=Fit}
\usepackage[T1]{fontenc}\usepackage{lmodern}

\usepackage{mathtools,amsfonts,amssymb}
\usepackage{old-arrows}
\usepackage{caption}
\usepackage{subcaption}
\usepackage{booktabs}
\usepackage{pifont}
\usepackage{placeins}

\definecolor{green}{rgb}{0,0.5,0}

\usepackage{tikz}
\usetikzlibrary{calc}
\usepackage{tkz-euclide}
\usepackage{calc}
\usepackage{tikz,pgfplots}
\pgfplotsset{compat=1.12}
\usetikzlibrary{arrows.meta,3d,matrix,positioning,decorations.markings,math}
\makeatletter
\tikzoption{canvas is xy plane at z}[]{%
  \def\tikz@plane@origin{\pgfpointxyz{0}{0}{#1}}%
  \def\tikz@plane@x{\pgfpointxyz{1}{0}{#1}}%
  \def\tikz@plane@y{\pgfpointxyz{0}{1}{#1}}%
  \tikz@canvas@is@plane
}
\makeatother

\makeatletter
\tikzoption{canvas is plane}[]{\@setOxy#1}
\def\@setOxy O(#1,#2,#3)x(#4,#5,#6)y(#7,#8,#9)%
  {\def\tikz@plane@origin{\pgfpointxyz{#1}{#2}{#3}}%
   \def\tikz@plane@x{\pgfpointxyz{#4}{#5}{#6}}%
   \def\tikz@plane@y{\pgfpointxyz{#7}{#8}{#9}}%
   \tikz@canvas@is@plane
  }
\makeatother

\usepackage{xfrac}

\newcommand{\pz}{\phantom{0}}

\newcommand{\vpad}{\vphantom{\bigg(}}

\newcommand{\maxpoly}{d}

\newcommand{\nelem}{N}

\newcommand{\nfunc}{{n_\func}}

\newcommand{\func}{f}

\newcommand{\npoints}{n}

\definecolor{orange}{rgb}{1,0.5,0}
\definecolor{green}{rgb}{0,0.5,0}
\definecolor{purple}{rgb}{0.5,0,0.5}
\newcommand{\reviewerOne}[1]{#1}
\newcommand{\reviewerTwo}[1]{#1}
\newcommand{\reviewerThree}[1]{#1}
\usepackage[margin=1in]{geometry}
\biboptions{sort&compress}
\begin{document}

\begin{center}
{\Large
A Note on the Convergence of Symmetric Triangle Quadrature Rules}
\\[2em]
Brian A.~Freno$^\text{a}$, Neil R.~Matula$^\text{a}$, Joseph E.~Bishop$^\text{a}$
\\[1em]
{\footnotesize\textit{%
$^\text{a}$Sandia National Laboratories, Albuquerque, NM 87185}}
\end{center}

%===============================================================================
\section{Introduction} %========================================================
%===============================================================================
\label{sec:introduction}

Symmetric quadrature rules for triangles are commonly used to efficiently integrate two-dimensional domains in finite-element-type problems.  Most quadrature rules have been formulated to exactly integrate polynomials.  Examples of such rules for triangles are provided in~\cite{lyness_1975,dunavant_1985,wandzura_2003,papanicolopulos_2015,papanicolopulos_2015_anc2}. Approaches to formulate rules for arbitrary function sequences have also been presented for one dimension~\cite{ma_1996}, as well as for triangles in two dimensions~\cite{freno_quad}.  Symmetric rules are generally preferable, as they yield straightforward mappings to the physical domain, avoiding the need to choose a particular vertex mapping.

While the development of polynomial rules focuses on the maximum degree a given number of points can exactly integrate, smooth integrands are generally not polynomials of finite degree.  Therefore, for such integrands, one needs to balance integration accuracy and computational cost.  A natural approach to this balance is to choose the number of points such that the convergence rate with respect to the mesh size $h$ matches that of the other properties of the scheme, such as the planar or curved triangles that approximate the geometry or the basis functions that approximate the solution.  In other words, the integration error should converge at least as quickly as the geometry- and solution-approximation errors.

In general, it is expected that a quadrature rule capable of integrating polynomials up to degree $\maxpoly$ yields an integration error that is $\mathcal{O}(h^p)$, where $p=\maxpoly+1$.  However, as we describe in this paper, for symmetric triangle quadrature rules, when $\maxpoly$ is even, $p=\maxpoly+2$; therefore, for a $p^\text{th}$-order-accurate quadrature rule, fewer quadrature points are necessary, reducing the time required for matrix assembly in finite-element-type problems.  This reduction in cost is modest for local differential operators that yield sparse matrices but appreciable for global integral operators that yield dense matrices.  An example of the latter occurs in the method-of-moments implementation of the electromagnetic field integral equations, where four-dimensional integrals are evaluated over two two-dimensional domains for every matrix element~\cite{peterson_1998}.

In this paper, we briefly summarize the details of symmetric triangle quadrature rules, discuss error implications for quadrature rules for one dimension and triangles, and we provide numerical examples that support our observation that polynomials that exactly integrate even maximum degrees converge faster than the conventional expectation for sequences of regular meshes.
%===============================================================================
\section{Symmetric Triangle Quadrature Rules} %=================================
%===============================================================================
\label{sec:quadrature}

A sequence of $\nfunc$ functions $\mathbf{f}(\mathbf{x})=\{f_1(\mathbf{x}),\hdots,f_\nfunc(\mathbf{x})\}$ is exactly integrated by an $\npoints$-point quadrature rule by taking a weighted combination of the function values at $\mathbf{x}_i$ for each of the $\npoints$ points:
\begin{align}
\int_A \mathbf{f}(\mathbf{x})dA = \sum_{i=1}^\npoints w_i \mathbf{f}(\mathbf{x}_i).
\label{eq:quadrature}
\end{align}
In one dimension, $\nfunc=2\npoints$ and, for polynomials, $\mathbf{f}(x)=\{1,\hdots,x^{2\npoints-1}\}$, such that the maximum polynomial degree integrated exactly is
\begin{align}
\maxpoly = 2\npoints -1.
\label{eq:d_1d}
\end{align}  
However, in two dimensions, determining an expression for $\nfunc$ is less straightforward~\cite{xiao_2010}.

% Symmetric Rules for Triangles ------------------------------------------------

\begin{figure}
\centering
\begin{tikzpicture}

% Preliminary
\def\ts{5}; % Triangle scale factor

\coordinate (A) at (-0.5,{-sqrt(3)/6});
\coordinate (B) at ( 0.5,{-sqrt(3)/6});
\coordinate (C) at ( 0  , {sqrt(1/3)});
\coordinate (D) at ($0.5*(A)+0.5*(B)$);
\coordinate (E) at ($0.5*(B)+0.5*(C)$);
\coordinate (F) at ($0.5*(C)+0.5*(A)$);
\coordinate (O) at (0,0);

% First ----------------------------------------------
\def\tx{0}
\def\ty{{\ts*sqrt(3)/6}}

% Triangle
\coordinate (T1) at ($\ts*(A)+(\tx,\ty)$);
\coordinate (T2) at ($\ts*(B)+(\tx,\ty)$);
\coordinate (T3) at ($\ts*(C)+(\tx,\ty)$);
\coordinate (T4) at ($\ts*(D)+(\tx,\ty)$);
\draw[thick] (T1) -- (T2) -- (T3) -- cycle;

% Medians
\draw[dashed] ($\ts*(C)+(\tx,\ty)$) -- ($\ts*(D)+(\tx,\ty)$);
\draw[dashed] ($\ts*(A)+(\tx,\ty)$) -- ($\ts*(E)+(\tx,\ty)$);
\draw[dashed] ($\ts*(B)+(\tx,\ty)$) -- ($\ts*(F)+(\tx,\ty)$);

% Nodes
\draw[fill=black] ($\ts*(0,0)+(\tx,\ty)$) circle (.07);

\node[below = .1 of T4] {(a)};

% Second ----------------------------------------------
\def\tx{5.5}
\def\ty{{\ts*sqrt(3)/6}}

% Triangle
\coordinate (T1) at ($\ts*(A)+(\tx,\ty)$);
\coordinate (T2) at ($\ts*(B)+(\tx,\ty)$);
\coordinate (T3) at ($\ts*(C)+(\tx,\ty)$);
\coordinate (T4) at ($\ts*(D)+(\tx,\ty)$);
\draw[thick] (T1) -- (T2) -- (T3) -- cycle;

% Medians
\draw[dashed] ($\ts*(C)+(\tx,\ty)$) -- ($\ts*(D)+(\tx,\ty)$);
\draw[dashed] ($\ts*(A)+(\tx,\ty)$) -- ($\ts*(E)+(\tx,\ty)$);
\draw[dashed] ($\ts*(B)+(\tx,\ty)$) -- ($\ts*(F)+(\tx,\ty)$);

% Nodes
\draw[fill=black] ($\ts*(-0.25,-0.144337567297406)+(\tx,\ty)$) circle (.07);
\draw[fill=black] ($\ts*( 0.25,-0.144337567297406)+(\tx,\ty)$) circle (.07);
\draw[fill=black] ($\ts*( 0.00, 0.288675134594813)+(\tx,\ty)$) circle (.07);

\node[below = .1 of T4] {(b)};

% Third ----------------------------------------------

\def\tx{11}
\def\ty{{\ts*sqrt(3)/6}}

% Triangle
\coordinate (T1) at ($\ts*(A)+(\tx,\ty)$);
\coordinate (T2) at ($\ts*(B)+(\tx,\ty)$);
\coordinate (T3) at ($\ts*(C)+(\tx,\ty)$);
\coordinate (T4) at ($\ts*(D)+(\tx,\ty)$);
\draw[thick] (T1) -- (T2) -- (T3) -- cycle;
% Medians
\draw[dashed] ($\ts*(C)+(\tx,\ty)$) -- ($\ts*(D)+(\tx,\ty)$);
\draw[dashed] ($\ts*(A)+(\tx,\ty)$) -- ($\ts*(E)+(\tx,\ty)$);
\draw[dashed] ($\ts*(B)+(\tx,\ty)$) -- ($\ts*(F)+(\tx,\ty)$);

% Nodes
\draw[fill=black] ($\ts*( 0.061447179740077, 0.282059528176809)+(\tx,\ty)$) circle (.07);
\draw[fill=black] ($\ts*(-0.274994306650607,-0.087814945442589)+(\tx,\ty)$) circle (.07);
\draw[fill=black] ($\ts*( 0.213547126910531,-0.194244582734219)+(\tx,\ty)$) circle (.07);
\draw[fill=black] ($\ts*( 0.274994306650607,-0.087814945442589)+(\tx,\ty)$) circle (.07);
\draw[fill=black] ($\ts*(-0.061447179740077, 0.282059528176809)+(\tx,\ty)$) circle (.07);
\draw[fill=black] ($\ts*(-0.213547126910531,-0.194244582734219)+(\tx,\ty)$) circle (.07);

\node[below = .1 of T4] {(c)};

\end{tikzpicture}
\caption{Examples of the (a) type-0, (b) type-1, and (c) type-2 orbits used to construct symmetric triangle quadrature rules.}
\label{fig:three_orbits}
\end{figure}
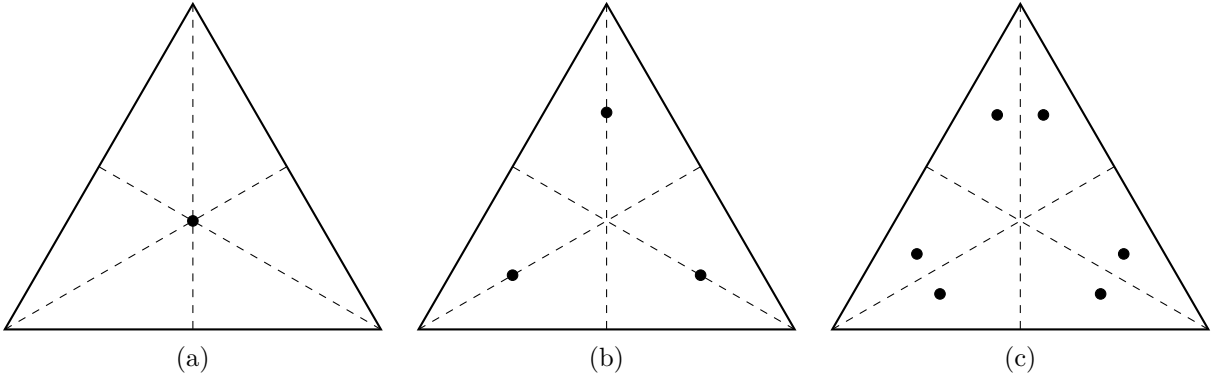

For triangles, \reviewerTwo{quadrature rules that are designed to be symmetric under the dihedral group $D_3$ are invariant under the three rotations and three reflections of an equilateral triangle, which can be isoparametrically mapped to an arbitrary triangle}.  Symmetric triangle rules are constructed using a combination of the three types of orbits~\cite{lyness_1975,dunavant_1985,wandzura_2003,papanicolopulos_2015}, examples of which are shown in Figure~\ref{fig:three_orbits}.

The type-0 orbit includes the centroid, which is $\left(\frac{1}{3},\frac{1}{3},\frac{1}{3}\right)$ in barycentric coordinates.  The type-1 orbit includes three points, each on a median, such that the coordinates are the three unique permutations of $\left(\lambda_1,\frac{1-\lambda_1}{2},\frac{1-\lambda_1}{2}\right)$.  Finally, the type-2 orbit includes six points, not on the medians; the coordinates are the six unique permutations of $(\lambda_1,\lambda_2,1-\lambda_1-\lambda_2)$.
The number of points $\npoints$ can be expressed as
\begin{align*}
\npoints = n_0 + 3n_1 + 6n_2,
\end{align*}
where $n_j$ is the number of type-$j$ orbits.   For each orbit, the one, three, or six points have the same weight.

% Polynomial Integration %------------------------------------------------------

Polynomial rules capable of integrating polynomials up to degree $\maxpoly$ can exactly integrate linear combinations of the monomials $x^p y^q$, where $0\le p\le \maxpoly$, $0\le q\le \maxpoly$, $0\le p+q\le \maxpoly$, totaling $\nfunc=(\maxpoly+1)(\maxpoly+2)/2$ monomials.  This definition yields more equations than unknowns; however, the equations are linearly dependent.

%===============================================================================
\section{Error Analysis} %======================================================
%===============================================================================
\label{sec:error}

We are interested in determining the convergence rate of the integration error 
\begin{align}
|e| \le C h^p
\label{eq:e}
\end{align}
over a discretized domain.  In~\eqref{eq:e}, $e=I-\tilde{I}$ is the integration error, $h$ is a one-dimensional measure of the mesh size, and $p$ is the rate of convergence.  $I = \int_A f(\mathbf{x})dA$ denotes the exact evaluation of the integral over the domain. $\tilde{I}$ denotes the quadrature evaluation of the integral over the domain by subdividing the domain into $\nelem$ elements and evaluating the integral over each element using quadrature:
\begin{align*}
\tilde{I} = \sum_{j=1}^{\nelem} \tilde{I}_j, \qquad \tilde{I}_j = \sum_{i=1}^\npoints w_{j_i} f(\mathbf{x}_{j_i}).
\end{align*}

%===============================================================================
\subsection{One-Dimensional Rules}
\label{sec:1d}

For one-dimensional polynomial quadrature rules, the integration error for a sufficiently smooth function $f$ over the domain $x\in[a,\,b]$ is given by \cite[Chap.~5]{kahaner_1989}
\begin{align}
e =  \int_a^b f(x) dx - \sum_{i=1}^\npoints w_i f(x_i) = \frac{(b-a)^{2\npoints+1}(\npoints!)^4}{(2\npoints+1)[(2\npoints)!]^3}f^{2\npoints}(\xi),
\label{eq:rn_1d}
\end{align}
where $\xi\in(a,b)$.  If the domain is subdivided into ${\nelem}_1$ elements of size $h_1=(b-a)/{\nelem}_1=h$, where quadrature is used to evaluate the integral over each element, the integration error over the domain is bounded by 
\begin{align}
|e_1| = \biggl|\sum_{j=1}^{{\nelem}_1} \frac{h^{2\npoints+1}(\npoints!)^4}{(2\npoints+1)[(2\npoints)!]^3}f^{2\npoints}(\xi_j)\biggr| \le C_1 h^{2\npoints},
\label{eq:e1}
\end{align}
where
\begin{align*}
C_1 = \frac{(b-a)(\npoints!)^4}{(2\npoints+1)[(2\npoints)!]^3}\max_{\xi\in(a,b)} f^{2\npoints}(\xi).
\end{align*}
If the domain is subdivided into ${\nelem}_2$ elements, where ${\nelem}_2=q {\nelem}_1$ and $h_2=h_1/q$, the integration error over the domain is bounded by 
\begin{align}
|e_2| = \biggl|\sum_{j=1}^{{\nelem}_2} \frac{h^{2\npoints+1}(\npoints!)^4}{q^{2\npoints+1}(2\npoints+1)[(2\npoints)!]^3}f^{2\npoints}(\xi_j)\biggr| \le \frac{C_1}{q^{2\npoints}} h^{2\npoints}.
\label{eq:e2}
\end{align}
From~\eqref{eq:e}, \eqref{eq:e1}, and \eqref{eq:e2}, as $h$ approaches zero,
\begin{align}
p = \frac{\log |e_1|/|e_2|}{\log q} \approx {2\npoints}.
\label{eq:p_1d}
\end{align}
From~\eqref{eq:d_1d} and~\eqref{eq:p_1d}, $p=\maxpoly+1$, as mentioned in the introduction.

%===============================================================================
\subsection{Triangular Rules}

\newcommand{\oidx}{{m}}%{{d'}}
\newcommand{\iidx}{{k}}%{{d''}}

For two dimensions, extensions of one-dimensional approaches are not straightforward, including the ability to derive an expression similar to~\eqref{eq:rn_1d}; however, in this work, we are interested in $p$ in~\eqref{eq:e}, not $C$.  

We consider the integral of $F(\mathbf{x})$ over a triangular element $\mathcal{T}$ of size $h$, which can be expressed in barycentric coordinates $(\lambda_1,\,\lambda_2,\,\lambda_3)$: 
\begin{align*}
\int_\mathcal{T} F(\mathbf{x}) dA = \frac{A}{A_\text{ref}}\int_0^1 \int_0^{1-\lambda_1} f(\boldsymbol{\lambda})d\lambda_2 d\lambda_1,  
\end{align*}
where $A$ is the area of the physical triangle and $A_\text{ref}=1/2$ is the area of the reference triangle.  Since $\lambda_3=1-\lambda_1-\lambda_2$, we only consider the independent coordinates $\boldsymbol{\lambda}=(\lambda_1,\,\lambda_2)$.  Finally, $f(\boldsymbol{\lambda})=F(\mathbf{x}(\boldsymbol{\lambda}))$.
The integration error over this triangular element is given by
\begin{align*}
e_h = \frac{A}{A_\text{ref}}\int_0^1 \int_0^{1-\lambda_1} f(\boldsymbol{\lambda})d\lambda_2 d\lambda_1 - \frac{A}{A_\text{ref}}\sum_{i=1}^\npoints w_{i} f(\boldsymbol{\lambda}_i),
\end{align*}
where the quadrature points $\boldsymbol{\lambda}_i$ and weights $w_{i}$ are expressed in barycentric coordinates.

A Taylor series expansion of $f(\lambda_1,\lambda_2)$ about the center of the reference triangle $\boldsymbol{\lambda}_c=(\lambda_c,\lambda_c)=(1/3,1/3)$ yields
\begin{align}
f(\lambda_1,\lambda_2) = \sum_{\oidx=0}^\infty \sum_{\iidx=0}^{\oidx} \frac{1}{(\oidx-\iidx)!(\iidx)!}\frac{\partial^{\oidx} f}{\partial \lambda_1^{\oidx-\iidx}\partial \lambda_2^{\iidx}}(\lambda_1-\lambda_c)^{\oidx-\iidx}(\lambda_2-\lambda_c)^{\iidx},
\label{eq:taylor}
\end{align}
where $\partial^{\oidx} f/(\partial \lambda_1^{\oidx-\iidx}\partial \lambda_2^{\iidx})$ is evaluated at $\boldsymbol{\lambda}_c$.
For an $\npoints$-point rule that can exactly integrate polynomials up to degree $\maxpoly$, the remainder of $f$ that cannot be integrated exactly is given by
\begin{align*}
r(\lambda_1,\lambda_2) 
&{}= \sum_{\oidx=\maxpoly+1}^\infty \sum_{\iidx=0}^{\oidx} \frac{1}{(\oidx-\iidx)!(\iidx)!}\frac{\partial^{\oidx} f}{\partial \lambda_1^{\oidx-\iidx}\partial \lambda_2^{\iidx}}(\lambda_1-\lambda_c)^{\oidx-\iidx}(\lambda_2-\lambda_c)^{\iidx}
\\[1em]
&{}= \sum_{\iidx=0}^{\maxpoly+1} \frac{1}{(\maxpoly+1-\iidx)!(\iidx)!}\frac{\partial^{\maxpoly+1} f}{\partial \lambda_1^{\maxpoly+1-\iidx}\partial \lambda_2^{\iidx}}(\lambda_1-\lambda_c)^{\maxpoly+1-\iidx}(\lambda_2-\lambda_c)^{\iidx}+\mathcal{O}(h^{\maxpoly+2}).
%\label{eq:taylor}
\end{align*}
%
%where $(\lambda_1-\lambda_c)$ and $(\lambda_2-\lambda_c)$ are each $\mathcal{O}(h)$ in physical coordinates.
%where $\hat{h}<1$, $\mathcal{O}(\hat{h})$ in barycentric coordinates is $\mathcal{O}(h)$ in physical coordinates, and $(\lambda_1-\lambda_c)$ and $(\lambda_2-\lambda_c)$ are each $\mathcal{O}(\hat{h})$.

\begin{figure}
\centering
\begin{tikzpicture}

\definecolor{lightgreen}{RGB}{173.4,217,170.6}
\definecolor{darkgreen} {RGB}{ 51  ,160, 44}
\definecolor{darkblue}  {RGB}{ 31  ,120,180}
\definecolor{lightblue} {RGB}{165.4,201,225}
\definecolor{darkred}   {RGB}{227  ,26,28}
\definecolor{lightred}  {RGB}{243.8,163.4,164.2}
\definecolor{darkorange}{RGB}{255  ,127  ,0}

% Preliminary
\def\ts{5}; % Triangle scale factor

\def\alp{0.7}
\def\bet{0.2}
\def\gam{0.4}

\coordinate (A) at (0,0); % Bottom Left
\coordinate (B) at (1,0); % Bottom Right
\coordinate (C) at (0,1); % Top Left

\coordinate (D) at (${(1-\alp)}*(A)+\alp*(B)$); % Bottom Edge
\coordinate (E) at (${(1-\gam)}*(C)+\gam*(B)$); % Diagonal Edge
\coordinate (F) at (${(1-\bet)}*(A)+\bet*(C)$); % Left Edge

\coordinate (G) at ($1/3*(B)+1/3*(E)+1/3*(D)$); % Center 1
\coordinate (H) at ($1/3*(E)+1/3*(C)+1/3*(F)$); % Center 2
\coordinate (I) at ($1/3*(D)+1/3*(F)+1/3*(A)$); % Center 3
\coordinate (J) at ($1/3*(F)+1/3*(D)+1/3*(E)$); % Center 4

\coordinate (O) at (0,0);

% First ----------------------------------------------
\def\tx{0}
\def\ty{{\ts*sqrt(3)/6}}

% Triangle
\coordinate (T1) at ($\ts*(A)+(\tx,\ty)$);
\coordinate (T2) at ($\ts*(B)+(\tx,\ty)$);
\coordinate (T3) at ($\ts*(C)+(\tx,\ty)$);
\coordinate (T4) at ($\ts*(D)+(\tx,\ty)$);
\coordinate (T5) at ($\ts*(E)+(\tx,\ty)$);
\coordinate (T6) at ($\ts*(F)+(\tx,\ty)$);
\coordinate (T7) at ($\ts*(G)+(\tx,\ty)$);
\coordinate (T8) at ($\ts*(H)+(\tx,\ty)$);
\coordinate (T9) at ($\ts*(I)+(\tx,\ty)$);
\coordinate (T0) at ($\ts*(J)+(\tx,\ty)$);

\node[above = 0 of T3] {$(0,1,0)$};
\node[anchor=north] at (T2) {$(1,0,0)$};
\node[anchor=north] at (T1) {$(0,0,1)$};

% Subtriangles
\draw[dashed] ($\ts*(D)+(\tx,\ty)$) -- ($\ts*(E)+(\tx,\ty)$);
\draw[dashed] ($\ts*(E)+(\tx,\ty)$) -- ($\ts*(F)+(\tx,\ty)$);
\draw[dashed] ($\ts*(F)+(\tx,\ty)$) -- ($\ts*(D)+(\tx,\ty)$);

% Nodes
%\draw[draw=black,fill=black] ($(T1)$) circle (.07);
%\draw[draw=black,fill=black] ($(T2)$) circle (.07);
%\draw[draw=black,fill=black] ($(T3)$) circle (.07);
%\draw[draw=black,fill=black] ($(T4)$) circle (.07);
%\draw[draw=black,fill=black] ($(T5)$) circle (.07);
%\draw[draw=black,fill=black] ($(T6)$) circle (.07);

\node[anchor=north] at (T4) {$(\alpha,0,1-\alpha)$};
\node[anchor=south west] at (T5) {$(\gamma,1-\gamma,0)$};
\node[anchor=east] at (T6) {$(0,\beta,1-\beta)$};

\node[anchor=center] at (T7) {$(1)$};
\node[anchor=center] at (T8) {$(2)$};
\node[anchor=center] at (T9) {$(3)$};
\node[anchor=center] at (T0) {$(4)$};

\draw[thick] (T1) -- (T2) -- (T3) -- cycle;
\end{tikzpicture}
\caption{Subdivision of the reference triangle in barycentric coordinates.}
\label{fig:subtriangles}
\end{figure}
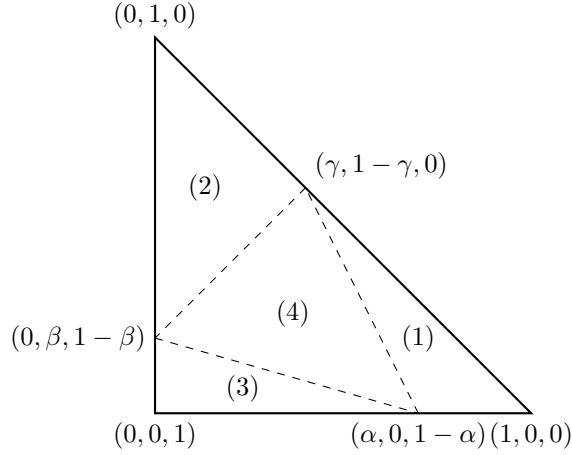

The leading term of the integration error is proportional to the integral of the leading term of the remainder:
\begin{align*}
e_h &{}= C 
\sum_{\iidx=0}^{\maxpoly+1} \frac{A}{(\maxpoly+1-\iidx)!(\iidx)!}\frac{\partial^{\maxpoly+1} f}{\partial \lambda_1^{\maxpoly+1-\iidx}\partial \lambda_2^{\iidx}}\int_0^1 \int_0^{1-\lambda_1} (\lambda_1-\lambda_c)^{\maxpoly+1-\iidx}(\lambda_2-\lambda_c)^{\iidx} d\lambda_2 d\lambda_1+\mathcal{O}(h^{\maxpoly+4})
\\
&{}= e_1 + \mathcal{O}(h^{\maxpoly+4}),
%\nonumber
%\\
%&{}= 
%C 
%\sum_{\iidx=0}^{d+1} B_{d,\iidx}\frac{\partial^{d+1} f}{\partial \lambda_1^{d+1-\iidx}\partial \lambda_2^{\iidx}}+\mathcal{O}(\hat{h}^{d+2})
%,
%\label{eq:error_h}
\end{align*}
where 
\begin{align}
e_1 = C 
\sum_{\iidx=0}^{\maxpoly+1} \frac{A}{(\maxpoly+1-\iidx)!(\iidx)!}\frac{\partial^{\maxpoly+1} f}{\partial \lambda_1^{\maxpoly+1-\iidx}\partial \lambda_2^{\iidx}}\int_0^1 \int_0^{1-\lambda_1} (\lambda_1-\lambda_c)^{\maxpoly+1-\iidx}(\lambda_2-\lambda_c)^{\iidx} d\lambda_2 d\lambda_1,
\label{eq:e_1}
\end{align}
and $A_\text{ref}$ is incorporated into $C$.

Next, we subdivide the reference triangle into four subtriangles $(i)$ using vertices positioned along the edges: $(\alpha,0)$, $(\gamma,1-\gamma)$, and $(0,\beta)$, as depicted in Figure~\ref{fig:subtriangles}.  Each subtriangle has its own barycentric coordinates $\boldsymbol{\lambda}^{(i)}$.  The transformations are
\begin{alignat*}{7}
\lambda_1 &{}= \lambda_1^{(1)} + \alpha(1-\lambda_1^{(1)}-\lambda_2^{(1)}) + \gamma\lambda_2^{(1)},\qquad &
\lambda_2 &{}= (1-\gamma)\lambda_2^{(1)} , \\
\lambda_1 &{}= \gamma\lambda_1^{(2)}, \qquad &
\lambda_2 &{}= (1-\gamma)\lambda_1^{(2)} + \beta(1-\lambda_1^{(2)}-\lambda_2^{(2)}) +\lambda_2^{(2)}, \\
\lambda_1 &{}= \alpha\lambda_1^{(3)}, \qquad &
\lambda_2 &{}= \beta\lambda_2^{(3)}, \\
\lambda_1 &{}= \gamma(1-\lambda_1^{(4)}-\lambda_2^{(4)}) +\alpha\lambda_2^{(4)}, \qquad &
\lambda_2 &{}= (1-\gamma)(1-\lambda_1^{(4)}-\lambda_2^{(4)}) +\beta\lambda_1^{(4)}.
\end{alignat*}

To combine and relate the integration errors over the subtriangles to the reference triangle, it is necessary to transform the derivatives and integrals in~\eqref{eq:e_1}.
Let $\alpha$, $\beta$, and $\gamma$ be expressed as perturbations about the edge midpoints:
\begin{align*}
\alpha = \frac{1}{2} + a\delta, \qquad
\beta = \frac{1}{2} +  b\delta, \qquad
\gamma = \frac{1}{2} + c\delta, \qquad
%\label{eq:perturbations}
\end{align*}
where $\delta\ll 1$ to yield a regular mesh refinement.
The determinants of the Jacobians $\mathbf{J}^{(i)} = \partial \boldsymbol{\lambda} / \partial \boldsymbol{\lambda}^{(i)}$ are 
%
%\begin{align*}
%\det \mathbf{J}^{(1)} = (1-\alpha)(1-\gamma), \qquad
%\det \mathbf{J}^{(2)} = \gamma(1-\beta),      \qquad
%\det \mathbf{J}^{(3)} = \alpha\beta,      \qquad
%\det \mathbf{J}^{(4)} = \beta\gamma + \alpha(1-\beta-\gamma),      
%\end{align*}
%
\begin{alignat*}{7}
\det \mathbf{J}^{(1)} &{}= (1-\alpha)(1-\gamma)                &&{}= \frac{1}{4}-\frac{1}{2}(a + c)\delta + a c\delta^2 &&{}= \frac{1}{4} + \mathcal{O}(\delta), \\
\det \mathbf{J}^{(2)} &{}= \gamma(1-\beta)                     &&{}= \frac{1}{4}-\frac{1}{2}(b - c)\delta - b c\delta^2 &&{}= \frac{1}{4} + \mathcal{O}(\delta), \\
\det \mathbf{J}^{(3)} &{}= \alpha\beta                         &&{}= \frac{1}{4}+\frac{1}{2}(a + b)\delta + a b\delta^2 &&{}= \frac{1}{4} + \mathcal{O}(\delta), \\
\det \mathbf{J}^{(4)} &{}= \beta\gamma + \alpha(1-\beta-\gamma)&&{}= \frac{1}{4}- (a b + a c - b c)\delta^2             &&{}= \frac{1}{4} + \mathcal{O}(\delta^2).
\end{alignat*}
%
%or, from~\eqref{eq:perturbations}, $\det \mathbf{J}^{(i)} = 1/4 + \mathcal{O}(\hat{h})$.
%
The derivatives are transformed according to
\begin{multline*}
\frac{\partial^{\maxpoly+1} f}{\partial {\lambda_1^{(i)}}^{\maxpoly+1-{\iidx}}\partial {\lambda_2^{(i)}}^{{\iidx}}} 
\\
=
\sum_{m=0}^{\maxpoly+1-{\iidx}} \sum_{j=0}^{{\iidx}} \binom{\maxpoly+1-{\iidx}}{m}\binom{{\iidx}}{j}
\biggl(\frac{\partial \lambda_1}{\partial \lambda_1^{(i)}}\biggr)^{\maxpoly+1-{\iidx}-m}
\biggl(\frac{\partial \lambda_2}{\partial \lambda_1^{(i)}}\biggr)^m
\biggl(\frac{\partial \lambda_1}{\partial \lambda_2^{(i)}}\biggr)^{{\iidx}-j}
\biggl(\frac{\partial \lambda_2}{\partial \lambda_2^{(i)}}\biggr)^{j}
\frac{\partial^{\maxpoly+1} f}{\partial \lambda_1^{\maxpoly+1-j-m}\partial \lambda_2^{j+m}},
\end{multline*}
%========================================
such that
\begin{align*}
\frac{\partial^{\maxpoly+1} f}{\partial {\lambda_1^{(1)}}^{\maxpoly+1-{\iidx}}\partial {\lambda_2^{(1)}}^{{\iidx}}}
&{}=
\sum_{j=0}^\iidx \binom{\iidx}{j} \biggl(\frac{1}{2} - a\delta\biggr)^{\maxpoly+1-\iidx}\biggl(\frac{1}{2} - c\delta\biggr)^{j}(c\delta-a\delta)^{\iidx-j}
\frac{\partial^{\maxpoly+1} f}{\partial \lambda_1^{\maxpoly+1-j}\partial \lambda_2^{j}}
\\
&{}= \frac{1}{2^{\maxpoly+1}}\frac{\partial^{\maxpoly+1} f}{\partial \lambda_1^{\maxpoly+1-{\iidx}}\partial \lambda_2^{{\iidx}}} + \mathcal{O}(\delta h^{\maxpoly+1}),
\\[1em]
\frac{\partial^{\maxpoly+1} f}{\partial {\lambda_1^{(2)}}^{\maxpoly+1-{\iidx}}\partial {\lambda_2^{(2)}}^{{\iidx}}}
&{}=
\sum_{m=0}^{\maxpoly+1-\iidx} \binom{\maxpoly+1-\iidx}{m} \biggl(\frac{1}{2}+c\delta\biggr)^{\maxpoly+1-\iidx-m}\biggl(\frac{1}{2}-b\delta\biggr)^\iidx(-b\delta-c\delta)^m \frac{\partial^{\maxpoly+1} f}{\partial \lambda_1^{\maxpoly+1-\iidx-m}\partial \lambda_2^{\iidx+m}}
\\
&{}= \frac{1}{2^{\maxpoly+1}}\frac{\partial^{\maxpoly+1} f}{\partial \lambda_1^{\maxpoly+1-{\iidx}}\partial \lambda_2^{{\iidx}}} + \mathcal{O}(\delta h^{\maxpoly+1}),
\\[1em]
\frac{\partial^{\maxpoly+1} f}{\partial {\lambda_1^{(3)}}^{\maxpoly+1-{\iidx}}\partial {\lambda_2^{(3)}}^{{\iidx}}}
&{}=
\biggl(\frac{1}{2}+a\delta\biggr)^{\maxpoly+1-\iidx}\biggl(\frac{1}{2}+b\delta\biggr)^\iidx \frac{\partial^{\maxpoly+1} f}{\partial \lambda_1^{\maxpoly+1-\iidx}\partial \lambda_2^{\iidx}}
\\
&{}= \frac{1}{2^{\maxpoly+1}}\frac{\partial^{\maxpoly+1} f}{\partial \lambda_1^{\maxpoly+1-{\iidx}}\partial \lambda_2^{{\iidx}}} + \mathcal{O}(\delta h^{\maxpoly+1}),
\\[1em]
\frac{\partial^{\maxpoly+1} f}{\partial {\lambda_1^{(4)}}^{\maxpoly+1-{\iidx}}\partial {\lambda_2^{(4)}}^{{\iidx}}}
&{}=
\sum_{m=0}^{\maxpoly+1-\iidx}\sum_{j=0}^\iidx  \binom{\maxpoly+1-\iidx}{m}\binom{\iidx}{j} \biggl(-\frac{1}{2}-c\delta\biggr)^{\maxpoly+1-\iidx-m}\biggl(-\frac{1}{2}+c\delta\biggr)^j 
\\ & \qquad\qquad \times
(a\delta-c\delta)^{\iidx-j}(b\delta+c\delta)^m\frac{\partial^{\maxpoly+1} f}{\partial \lambda_1^{\maxpoly+1-j-m}\partial \lambda_2^{j+m}}
\\
&{}= \frac{(-1)^{\maxpoly+1}}{2^{\maxpoly+1}}
\frac{\partial^{\maxpoly+1} f}{\partial \lambda_1^{\maxpoly+1-{\iidx}}\partial \lambda_2^{{\iidx}}} + \mathcal{O}(\delta h^{\maxpoly+1}).
\end{align*}
Additionally, to combine the integration errors from the subtriangles, we take zeroth-order Taylor series expansions about the reference triangle center $\boldsymbol{\lambda}_c$ to approximate the derivatives at the subtriangle centers $\boldsymbol{\lambda}_c^{(i)}$:
\begin{align*}
\left.\frac{\partial^{\maxpoly+1} f}{\partial {\lambda_1^{(i)}}^{\maxpoly+1-{\iidx}}\partial {\lambda_2^{(i)}}^{{\iidx}}}\right|_{\boldsymbol{\lambda}_c^{(i)}}
&{}=
B_{\maxpoly+1}^{(i)}
\left.\frac{\partial^{\maxpoly+1} f}{\partial \lambda_1^{\maxpoly+1-{\iidx}}\partial \lambda_2^{{\iidx}}}\right|_{\boldsymbol{\lambda_c}}  + \mathcal{O}(\delta h^{\maxpoly+1}) + \mathcal{O}(h^{\maxpoly+2}),
\end{align*}
where
\begin{align*}
B_{\maxpoly+1}^{(1)} = \frac{1}{2^{\maxpoly+1}},\qquad
B_{\maxpoly+1}^{(2)} = \frac{1}{2^{\maxpoly+1}},\qquad
B_{\maxpoly+1}^{(3)} = \frac{1}{2^{\maxpoly+1}},\qquad
B_{\maxpoly+1}^{(4)} = \frac{(-1)^{\maxpoly+1}}{2^{\maxpoly+1}}.
\end{align*}

For each subtriangle, the error can be expressed in terms of $e_1$~\eqref{eq:e_1}:
%
%\begin{align}
%e_{h/2}^{(i)} = e_1 \det \mathbf{J}^{(i)} \left.\Biggl(\displaystyle\frac{\partial^{\maxpoly+1} f}{\partial {\lambda_1^{(i)}}^{\maxpoly+1-{\iidx}}\partial {\lambda_2^{(i)}}^{{\iidx}}}\Biggr)
%\middle/
%\biggl(\frac{\partial^{\maxpoly+1} f}{\partial \lambda_1^{\maxpoly+1-{\iidx}}\partial \lambda_2^{{\iidx}}}\biggr)
%\right. + \mathcal{O}(h^{\maxpoly+4})
%,
%\label{eq:error_local}
%\end{align}
%
\begin{align}
e_{h/2}^{(i)} = e_1 \det \mathbf{J}^{(i)} B_{\maxpoly+1}^{(i)} +\mathcal{O}(\delta h^{\maxpoly+1}) + \mathcal{O}(h^{\maxpoly+2}) 
=
e_1 \det \mathbf{J}^{(i)} B_{\maxpoly+1}^{(i)} +\mathcal{O}(h^{\maxpoly+1})
,
\label{eq:error_local}
\end{align}
where $\det \mathbf{J}^{(i)}=A^{(i)}/A$.
The total numerical integration error from the subtriangles is
%
%\begin{align*}
%e_{h/2}^{\vphantom(} =\sum_{i=1}^4 e_{h/2}^{(i)} = \biggl(\frac{1}{2^{\maxpoly+3}}\bigl(3+(-1)^{\maxpoly+1}\bigr) + \mathcal{O}(\hat{h})\biggr) e_1 = e_2 + \mathcal{O}(h^{\maxpoly+4}).
%\end{align*}
%
\begin{align*}
e_{h/2}^{\vphantom(} =\sum_{i=1}^4 e_{h/2}^{(i)} = \frac{1}{2^{\maxpoly+3}}\bigl(3+(-1)^{\maxpoly+1}\bigr) e_1 + \mathcal{O}(h^{\maxpoly+1}) = e_2 + \mathcal{O}(h^{\maxpoly+1}).
\end{align*}
where
\begin{align}
\frac{e_2}{e_1} = \frac{1}{2^{\maxpoly+3}}\bigl(3+(-1)^{\maxpoly+1}\bigr)= \left\{\begin{array}{@{} c l @{}}2^{-(\maxpoly+1)} & \text{$\maxpoly$ is odd} \\ 2^{-(\maxpoly+2)} & \text{$\maxpoly$ is even}\end{array}\right..
\label{eq:e_2}
\end{align}
Inserting~\eqref{eq:e_1} and~\eqref{eq:e_2} into~\eqref{eq:p_1d} and noting that the subdivision corresponds to $q=2$ shows 
\begin{align}
p = \frac{\log |e_1|/|e_2|}{\log q} = \left\{\begin{array}{@{} c l @{}}\maxpoly+ 1 & \text{$\maxpoly$ is odd} \\ \maxpoly+ 2 & \text{$\maxpoly$ is even}\end{array}\right..
\label{eq:p_2d}
\end{align}
Equation~\eqref{eq:p_2d} indicates that, when $\maxpoly$ is odd, $p=\maxpoly+1$, as seen in Section~\ref{sec:1d}.  However, when $\maxpoly$ is even, $p=\maxpoly+2$.  
Therefore, if a rule is selected based on a desired convergence rate $p$, selecting a rule that exactly integrates an even $\maxpoly$ requires fewer quadrature points and incurs less computational expense.

\reviewerTwo{%
It is important to note that the identical mixed derivatives and equal-magnitude coefficients in the leading terms of the subtriangle errors~\eqref{eq:error_local} are due to the $D_3$ symmetry of the quadrature rules.  In the absence of such symmetry, the uniform structure and error reduction would not be guaranteed.
Furthermore, had we summed the absolute values of $e_{h/2}^{(i)}$ instead of the signed values, we would have obtained the more pessimistic $p=\maxpoly+1$ when $\maxpoly$ is even.
}

%===============================================================================
\section{Numerical Examples} %==================================================
%===============================================================================
\label{sec:results}

\begin{figure}[p]
\centering
\includegraphics[scale=.28,clip=true,trim=1.5in 0in 1.5in 0in]{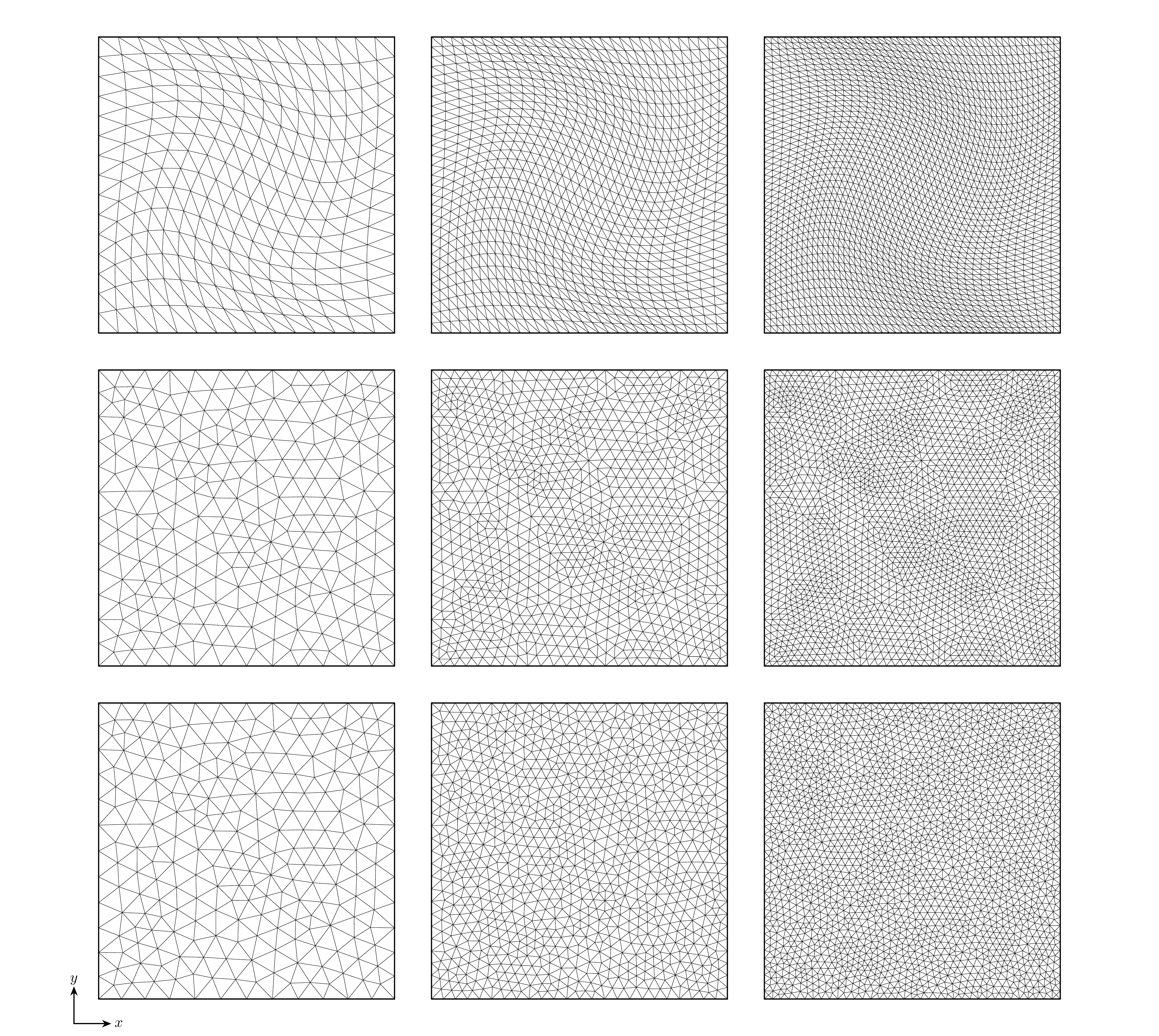}
\caption{The three coarsest meshes of each sequence: structured (top), nested unstructured (middle), and independent unstructured (bottom) meshes.}
\label{fig:meshes}
\end{figure}

\begin{table}[p]
\centering
\begin{tabular}{c c c c c c c c c c c c}
\toprule
$\maxpoly$ &  1 &  2 &  3 &  4 &  5 & \pz6 & \pz7 & \pz8 & \pz9 & 10 & 11 \\ \midrule
$\npoints$ &  1 &  3 &  4 &  6 &  7 &   12 &   13 &   16 &   19 & 25 & 27 \\
$p$        &  2 &  4 &  4 &  6 &  6 & \pz8 & \pz8 &   10 &   10 & 12 & 12 \\
\reviewerThree{$Q$}        & \reviewerThree{PI} & \reviewerThree{PI} & \reviewerThree{NI} & \reviewerThree{PI} & \reviewerThree{PI} &   \reviewerThree{PI} &   \reviewerThree{NI} &   \reviewerThree{PI} &   \reviewerThree{PI} & \reviewerThree{PI} & \reviewerThree{PO} \\
\bottomrule
\end{tabular}
\caption{Symmetric triangle quadrature properties: maximum polynomial degree $\maxpoly$, number of quadrature points $\npoints$, convergence rate $p$, \reviewerThree{and quality $Q$ (as described in Section~\ref{sec:results})}.}
\vskip-\dp\strutbox
\label{tab:quadrature}
\end{table}

In this section, we demonstrate the validity of~\eqref{eq:p_2d} for three sequences of ten meshes using quadruple precision, where the $k^\text{th}$ mesh consists of $\nelem=450 k^2$ triangles.
The first sequence consists of structured meshes, which are obtained by transforming uniform meshes using the transformation provided in~\cite{freno_2021}.  The second sequence consists of unstructured meshes that are nested refinements of the coarsest mesh, where each of the coarsest mesh triangles has been subdivided into $k^2$ similar subtriangles.  Finally, the third sequence consists of independently generated unstructured meshes.  The meshes in the first two sequences are systematically refined, whereas those in the third are not.
Figure~\ref{fig:meshes} shows the three coarsest meshes of each sequence.
  %The $k^\text{th}$ mesh consists of $\nelem=450 k^2$ triangles.  The structured meshes are obtained by transforming uniform meshes using the transformation provided in~\cite{freno_2021}, whereas the unstructured meshes are nested refinements of the coarsest mesh, where each of the coarsest mesh triangles has been subdivided into $k^2$ similar subtriangles.  The first three meshes of each type are shown in Figure~\ref{fig:meshes}.

\begin{figure}[!t]
\centering
\begin{subfigure}[b]{.49\textwidth}
\includegraphics[scale=.64,clip=true,trim=2.3in 0in 2.8in 0in]{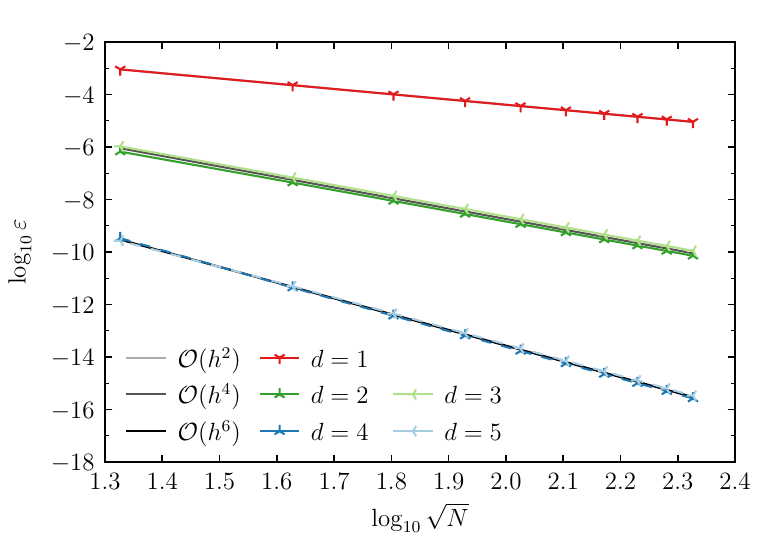}
\caption{Structured, $1\le\maxpoly\le5$\vpad}
\label{fig:str_1}
\end{subfigure}
\hspace{0.25em}
\begin{subfigure}[b]{.49\textwidth}
\includegraphics[scale=.64,clip=true,trim=2.3in 0in 2.8in 0in]{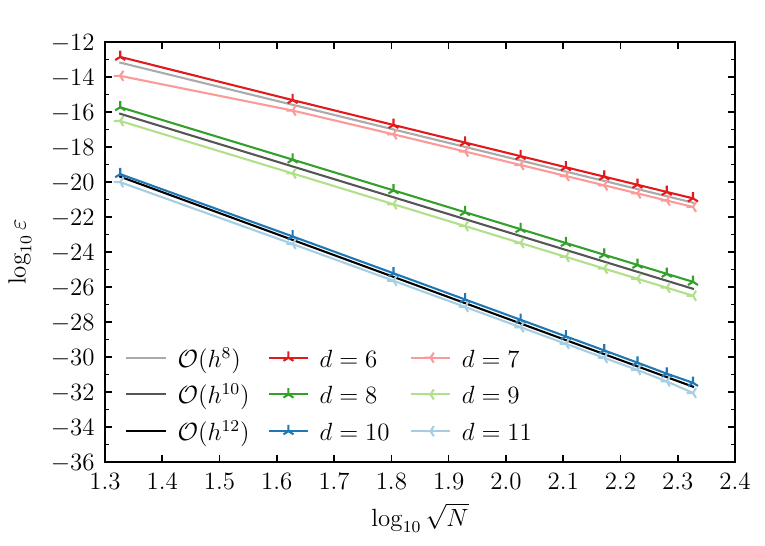}
\caption{Structured, $6\le\maxpoly\le11$\vpad}
\label{fig:str_2}
\end{subfigure}
\\
\begin{subfigure}[b]{.49\textwidth}
\includegraphics[scale=.64,clip=true,trim=2.3in 0in 2.8in 0in]{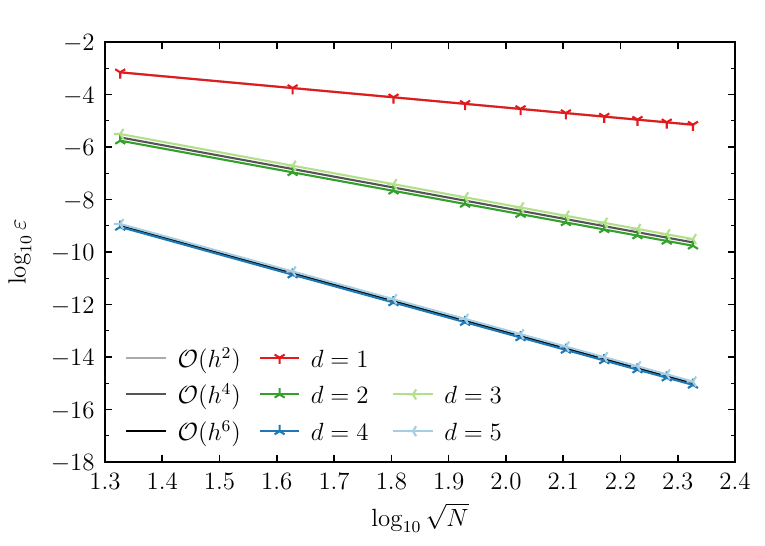}
\caption{Nested unstructured, $1\le\maxpoly\le5$\vpad}
\label{fig:uns_1}
\end{subfigure}
\hspace{0.25em}
\begin{subfigure}[b]{.49\textwidth}
\includegraphics[scale=.64,clip=true,trim=2.3in 0in 2.8in 0in]{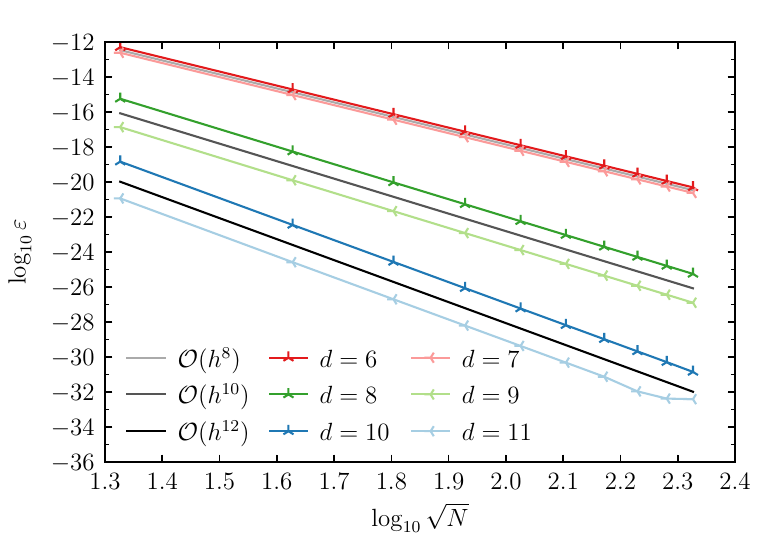}
\caption{Nested unstructured, $6\le\maxpoly\le11$\vpad}
\label{fig:uns_2}
\end{subfigure}
\\
\begin{subfigure}[b]{.49\textwidth}
\includegraphics[scale=.64,clip=true,trim=2.3in 0in 2.8in 0in]{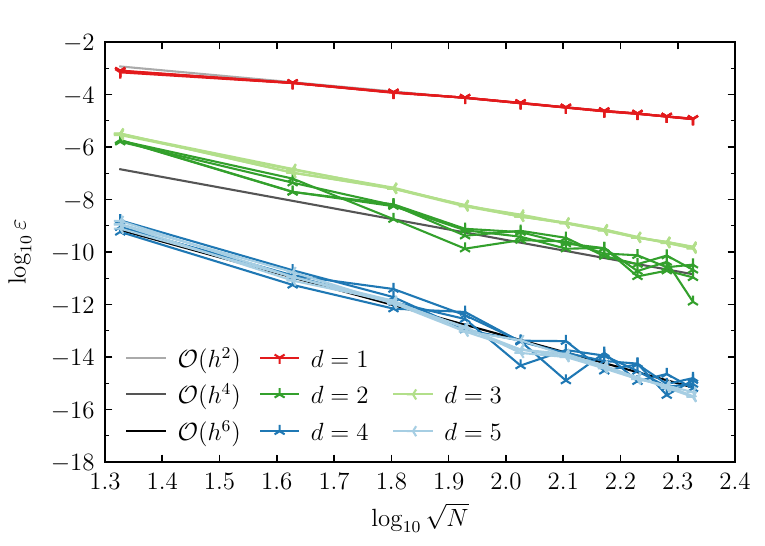}
\caption{Independent unstructured, $1\le\maxpoly\le5$\vpad}
\label{fig:uni_1}
\end{subfigure}
\hspace{0.25em}
\begin{subfigure}[b]{.49\textwidth}
\includegraphics[scale=.64,clip=true,trim=2.3in 0in 2.8in 0in]{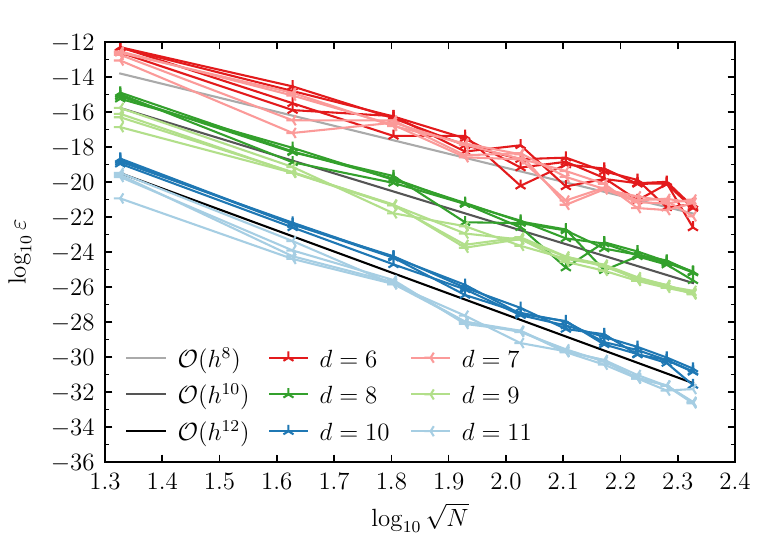}
\caption{Independent unstructured, $6\le\maxpoly\le11$\vpad}
\label{fig:uni_2}
\end{subfigure}
\caption{Integration error $\varepsilon=|e/I|$ \reviewerOne{with respect to the number of triangles $N$} for the first eleven quadrature rules for each mesh sequence.}
\vskip-\dp\strutbox
\label{fig:error}
\end{figure}

\begin{figure}[!t]
\centering
\begin{subfigure}[b]{.49\textwidth}
\includegraphics[scale=.64,clip=true,trim=2.3in 0in 2.8in 0in]{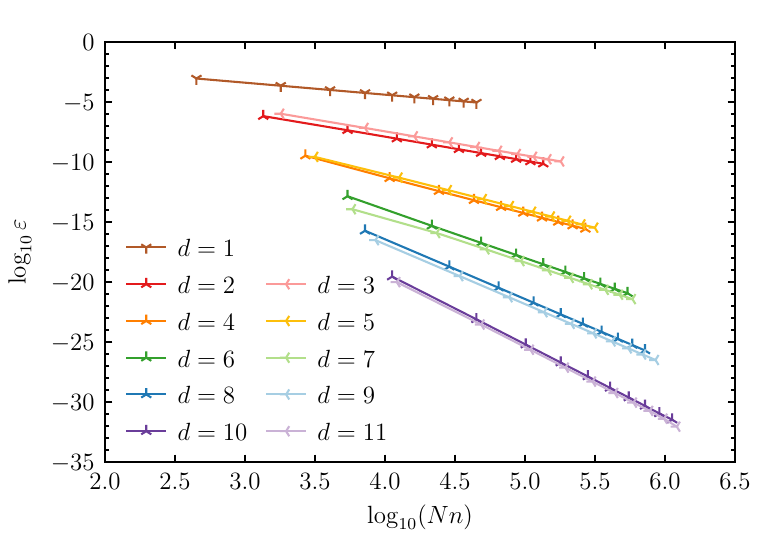}
\caption{Structured\vpad}
\label{fig:str_3_points}
\end{subfigure}
\hspace{0.25em}
\begin{subfigure}[b]{.49\textwidth}
\includegraphics[scale=.64,clip=true,trim=2.3in 0in 2.8in 0in]{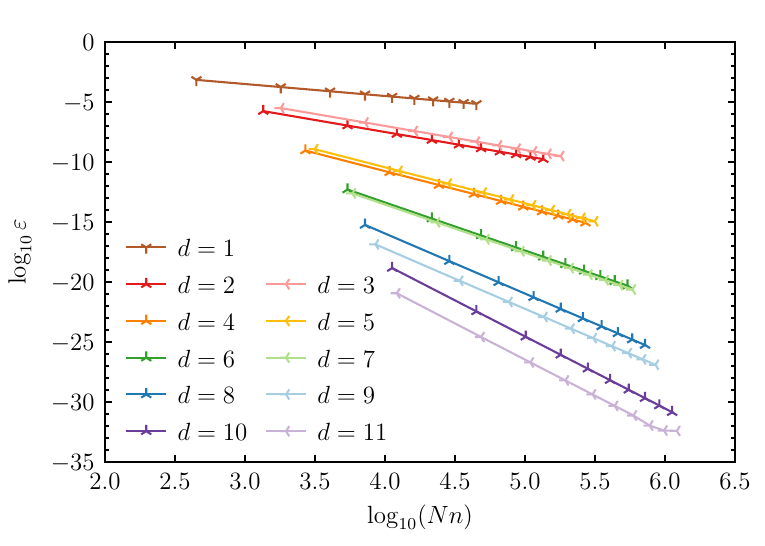}
\caption{Nested unstructured\vpad}
\label{fig:uns_3_points}
\end{subfigure}
\caption{ \reviewerOne{Integration error $\varepsilon=|e/I|$ with respect to the total number of quadrature evaluations $N n$ for the first eleven quadrature rules for each mesh sequence.}}
\vskip-\dp\strutbox
\label{fig:error_points}
\end{figure}

To observe large values of $p$, double precision is insufficient; the round-off error quickly exceeds the integration error.  Therefore, we use quadruple precision in this paper.  To achieve such precision for the quadrature rules, we numerically solve a nonlinear least-squares problem in Mathematica using the lower-precision representation of the rules provided in~\cite{dunavant_1985} as initial guesses.
Examples of structured meshes with double precision for the electromagnetic field integral equations are included in~\cite{freno_mfie_2022,freno_cfie_2023,freno_efie_slot_2024,freno_deep_slot_2025}.

We consider the integral over the unit square $(x,y)\in[0,1]\times[0,1]$ of 
\begin{align*}
f(x,y) = \exp(x+2 y) \bigl(1+\sin(3x+4y)+\sin(6x+5y)\bigr).
\end{align*}
This choice of function is arbitrary; however, it satisfies the criteria for a suitable function to demonstrate quadrature convergence: it is infinitely differentiable, there is no symmetry over the domain, it cannot be expressed as a finite-degree polynomial, and its integral over the domain can be computed analytically:
\begin{align*}
I ={}& \frac{1}{21460}\bigl(9317 \\
   &+e^{\phantom{1}} (-10730 + 1073\cos3 + 340\cos6\pz + 1073\sin3 + 560\sin6\pz) \\
   &+e^2             (-10730 + 1073\cos4 + 340\cos5\pz + 1073\sin4 + 560\sin5\pz) \\
   &-e^3             (-10730 + 1073\cos7 + 340\cos11   + 1073\sin7 + 560\sin11)\bigr).
\end{align*}

Figure~\ref{fig:error} shows the integration error $\varepsilon=|e/I|$ for each of the three mesh sequences arising from the first eleven quadrature rules \reviewerThree{from~\cite{dunavant_1985}, properties of which are listed in Table~\ref{tab:quadrature}.  The qualities $Q$ considered include positive weights with all points inside the triangle (PI), positive weights with some points outside the triangle (PO), and some negative weights with all points inside the triangle (NI).  These rules do not include the quality of some negative weights with some points outside the triangle (NO).  Provided the integrand over the triangle is sufficiently smooth and can be extrapolated beyond the triangle, these qualities will not affect the convergence rates.}
The integration errors clearly converge according to~\eqref{eq:p_2d} for the structured meshes in Figures~\ref{fig:str_1} and~\ref{fig:str_2} and the nested unstructured meshes in Figures~\ref{fig:uns_1} and~\ref{fig:uns_2}.  Unlike those sequences of meshes, which are systematically refined, we would generally not expect the independent unstructured meshes in the third sequence to converge as clearly.  For this sequence, we consider the four 90-degree rotations of each mesh.  Figures~\ref{fig:uni_1} and~\ref{fig:uni_2} show that the integration errors still converge according to~\eqref{eq:p_2d}.

\reviewerOne{%
Though the focus of this work is on convergence with respect to the number of triangles $N$, we show in Figure~\ref{fig:error_points} how the integration error $\varepsilon=|e/I|$ varies with respect to the total number of quadrature evaluations $N n$.  In Figure~\ref{fig:str_3_points}, for a fixed $p$, the two rule options perform similarly; however, the error for the rule with the larger number of points $\npoints$ that exactly integrates a higher degree $\maxpoly$ is lower for larger $p$.  In Figure~\ref{fig:uns_3_points}, this difference is more pronounced, especially for $\maxpoly=10$ and $\maxpoly=11$, but such observations may not hold for other integrands.  
Furthermore, in practice, the number of triangles is typically driven by other considerations, such as the desired geometry and solution approximations.
}

%===============================================================================
\section{Conclusions} %=========================================================
%===============================================================================
\label{sec:conclusions}

In this paper, we demonstrated that, for symmetric triangle quadrature rules that exactly integrate polynomials up to degree $\maxpoly$, the integration error is $\mathcal{O}(h^p)$, where $p=\maxpoly+1$ when $\maxpoly$ is odd but $p=\maxpoly+2$ when $\maxpoly$ is even.  Therefore, if a rule is selected based on a desired convergence rate, selecting a rule that exactly integrates an even degree requires fewer quadrature points and incurs less computational expense.
%===============================================================================
\section*{Acknowledgments} %====================================================
%===============================================================================
\label{sec:acknowledgments}

This article has been authored by employees of National Technology \& Engineering Solutions of Sandia, LLC under Contract No.~DE-NA0003525 with the U.S.~Department of Energy (DOE). The employees own all right, title, and interest in and to the article and are solely responsible for its contents. The United States Government retains and the publisher, by accepting the article for publication, acknowledges that the United States Government retains a non-exclusive, paid-up, irrevocable, world-wide license to publish or reproduce the published form of this article or allow others to do so, for United States Government purposes. The DOE will provide public access to these results of federally sponsored research in accordance with the DOE Public Access Plan \url{https://www.energy.gov/downloads/doe-public-access-plan}.

\addcontentsline{toc}{section}{\refname}
\bibliographystyle{elsarticle-num}
\bibliography{../quadrature_manuscript/quadrature.bib}
%\end{document}

%\input{notes.tex}

\end{document}